\documentclass[12pt,draft]{amsart}

\makeatletter

\theoremstyle{plain}    
\newtheorem{thm}{Theorem}[section]
\numberwithin{equation}{section} 
\numberwithin{figure}{section} 
\theoremstyle{plain}    
\newtheorem*{thm*}{Theorem} 
\theoremstyle{plain}    
\newtheorem{cor}[thm]{Corollary} 
\theoremstyle{plain}    
\newtheorem{lem}[thm]{Lemma} 
\theoremstyle{plain}    
\newtheorem{prop}[thm]{Proposition} 
\theoremstyle{definition}
\newtheorem{defn}[thm]{Definition}
\theoremstyle{remark}
\newtheorem{rem}[thm]{Remark}
\theoremstyle{remark}

\theoremstyle{remark}    

\theoremstyle{remark}    

\theoremstyle{definition}  

\theoremstyle{remark}
  \newtheorem*{acknowledgement*}{Acknowledgement} 

\theoremstyle{plain}    
\newtheorem{subthm}{Theorem}[subsection]
\theoremstyle{plain}    
\theoremstyle{plain}    
\theoremstyle{plain}    
\theoremstyle{definition}
\newtheorem{subdefn}[subthm]{Definition}
\theoremstyle{remark}

\theoremstyle{remark}    

\theoremstyle{remark}    

\theoremstyle{plain}    

\usepackage{amssymb} \usepackage{amscd}

\makeatother

\begin{document}

\title[Finite Entropy]{Finite Free Entropy and Free Group Factors}

\author{Nathanial P. Brown}

\address{Department of Mathematics, Penn State University, State
College, PA 16802}

\email{nbrown@math.psu.edu}

\thanks{Partially supported by DMS-0244807.}

\begin{abstract}
We show the existence of noncommutative random variables with finite
free entropy but which do not generate a free group
factor.  In particular, this gives an example of variables
$X_1,\ldots,X_n$ such that $\delta(X_1,\ldots,X_n) = n$ while
$W^*(X_1,\ldots,X_n) \ncong L({\mathbb F}_n)$.
\end{abstract}

\maketitle

\section{Introduction}

In the last few years we have seen deep analogies and exciting
interactions between the worlds of ergodic equivalence relations and
von Neumann algebras (cf.\ \cite{connes-dima}, \cite{gaboriau:betti},
\cite{gaboriau}, \cite{mineyev-dima}, \cite{monod-shalom},
\cite{ozawa:solid}, \cite{ozawa-popa}, \cite{popa:HT}). In particular,
there have been many parallels between the results obtained using
Voiculescu's theory of free entropy and free entropy dimension
(\cite{voiculescu:II}, \cite{kenley:packing}) and Gaboriau's work on
cost and $L^2$ Betti numbers of equivalence relations
\cite{gaboriau:betti}, \cite{gaboriau}.

One beautiful result in the realm of equivalence relations gives a
characterization of those which are ``treeable'' (i.e.\ a free product
of hyperfinite equivalence relations): they are precisely those
relations for which the cost can be attained \cite{gaboriau}. In other
words, an equivalence relation is a free product of hyperfinite
equivalence relations if and only if it has an ``optimal'' family of
generators, in that there is no family of generators with smaller
total ``support''.

The following question, pointed out to us by Dima Shlyakhtenko, is a
natural analogue of Gaboriau's characterization of treeable
equivalence relations.

\begin{center}
{\em If $\delta(X_1,\ldots,X_n) = n$ does it follow that
$W^*(X_1,\ldots,X_n)$ is a free group factor?}
\end{center}

Here $\delta$ denotes Voiculescu's free entropy dimension and 
$X_1,\ldots,X_n$ are self-adjoints in a tracial W$^*$-probability space.

To see that this is analogous to Gaboriau's characterization we first
note that the inequality $\delta(Y_1,\ldots,Y_m) \leq m$ always holds.
Hence, if it is true that $\delta$ is a W$^*$-invariant and
$Y_1,\ldots,Y_m$ is any other generating set for $W^*(X_1,\ldots,X_n)$
then the assumption $\delta(X_1,\ldots,X_n) = n$ would imply that $m
\geq n$.  Thus $X_1,\ldots,X_n$ would be a set of generators for
$W^*(X_1,\ldots,X_n)$ of minimal ``size'' and this minimum would be
attained -- in other words, the assumption $\delta(X_1,\ldots,X_n) =
n$ is a von Neumann algebraic analogue of the assumption that the cost
of a particular equivalence relation is attained.  So it is natural to
wonder if this assumption implies that $W^*(X_1,\ldots,X_n)$ is a free
product of hyperfinite von Neumann algebras.  If so, then it is
necessarily an (interpolated) free group factor (cf.\
\cite{voiculescu:circular}, \cite{dykema:freeproduct}) and one would
expect it to be $L({\mathbb F}_n)$.

It follows from work of Voiculescu (see \cite{voiculescu:survey} for a
nice survey) that if free entropy is finite (i.e.\
$\chi(X_1,\ldots,X_n) > -\infty$) then $\delta(X_1,\ldots,X_n) = n$
and hence the following question of Shlyakhtenko is weaker than the
one posed above. 

\begin{center}
{\em If $\chi (X_1,\ldots,X_n) > -\infty$ is $W^*(X_1,\ldots,X_n)$
necessarily a free group factor?}
\end{center}

In this paper we will give counterexamples to this question.  Of
course, they are also counterexamples to the first question as well.
More precisely, we will show that there are $n$-tuples
$X_1,\ldots,X_n$ such that $\chi (X_1,\ldots,X_n) > -\infty$ -- hence
$\delta (X_1,\ldots,X_n) = n$ -- but $W^*(X_1,\ldots,X_n)$ can't be
embed into a free group factor.  The reason no such embedding exists
is because we can arrange that $W^*(X_1,\ldots,X_n)$ has Haagerup
invariant strictly bigger than one and/or does not have Haagerup's
approximation property. Interestingly enough, however, Ozawa has
pointed out that our examples can be taken solid in the sense of
\cite{ozawa:solid}.  Of course, if one modifies the questions above by
further assuming that $W^*(X_1,\ldots,X_n)$ has the Haagerup
approximation property then they are both still open as Haagerup's
property and solidity are currently the only known obstructions to
embeddability into a free group factor.

In the next section we define all the concepts and review the known
results that we will need.  The third section contains the main
technical result which is inspired by (but does not actually use) the
theory of exact C$^*$-algebras.  The fourth section explains how to
combine the result from section three with other known results to get
counterexamples to the questions above.  In the final section of the
paper we observe that an affirmative answer to a considerable
weakening of the questions above would have important consequences both
inside and outside of free probability theory.

\begin{acknowledgement*}
Many thanks to my friends and colleagues Michael Anshelevich, Narutaka
Ozawa and Dima Shlyakhtenko for helpful remarks and comments.  I am
especially grateful to Dima for explaining the analogy with
equivalence relations and allowing me to plagiarize large parts of an
email he sent me in the introduction to this paper! Finally, I thank
Thomas Sch\"{u}cker and Antony Wassermann for the invitation to a
conference in Luminy where Dima informed me of the questions addressed
in this paper.
\end{acknowledgement*}

\section{Notation, Definitions and Known Results}

Throughout this paper we will follow standard notation in operator
algebra theory.  For a II$_1$-factor $M$, $L^2(M)$ will  denote the
GNS space of $M$ with respect to its unique tracial state.  An
inclusion $M \subset B(L^2(M))$ will always be assumed standard (i.e.\
arising from the GNS construction).  We let $M_*$ be the predual of
$M$ and $\|\cdot\|_2$ denote the 2-norm (i.e.\ $\|x\|_2 =
\sqrt{\tau(x^* x)}$, where $\tau$ is the trace on $M$). All von Neumann 
algebras are assumed to have separable predual. If $\Gamma$ is a
discrete (countable) group then we will denote by $L(\Gamma)$ the von
Neumann algebra generated by the left regular representation of
$\Gamma$.

If $R$ denotes the hyperfinite II$_1$-factor and
$\omega \in \beta({\mathbb N})\backslash {\mathbb N}$ is a free
ultrafilter then the ultraproduct $R^{\omega}$ is defined to be
$l^{\infty}(R) = \{ (x_n) : x_n \in R, \sup_n \| x_n \| < \infty \}$
modulo the ideal $I_{\omega} = \{ (x_n) : \lim_{n \to \omega} \| x_n
\|_2 = 0 \}$. As is well known, but not trivial, $R^{\omega}$ is a
II$_1$-factor with tracial state $\tau_{\omega} ((x_n)) = \lim_{n \to
\omega} \tau(x_n)$.

We will use the abbreviation u.c.p.\ for unital completely positive
maps.  The term completely bounded will be shortened to c.b.\ and the
c.b.\ norm of a map $\phi$ will be denoted $\|\phi\|_{cb}$ (see
\cite{paulsen:book} for more). An operator system is a norm closed,
self-adjoint subspace of a von Neumann algebra which contains the
unit.

\subsection{The Haagerup Approximation Property}

The Haagerup approximation property was formally introduced by Marie
Choda who was inspired by the fundamental paper
\cite{haagerup:freegroup} where Haagerup proved that free group
factors have this property.

\begin{subdefn} 
A finite von Neumann algebra $(M,\tau)$ is said to have the {\em
Haagerup approximation property} if there exists a sequence of u.c.p.\
$\tau$-preserving maps $\phi_n:M \to M$ such that (a) $\|\phi_n(x) -
x\|_2 \to 0$ for all $x \in M$ and (b) the map induced by $\phi_n$ on
$L^2(M,\tau)$ is a compact operator.
\end{subdefn}

The Cauchy-Schwartz inequality for u.c.p.\ maps implies that any trace
decreasing u.c.p.\ map $\phi:M\to M$ extends to a contraction on
$L^2(M,\tau)$ since $$\|\phi(x)\|_2^2 = \tau(\phi(x^*)\phi(x)) \leq
\tau(\phi(x^* x)) \leq \tau(x^* x) = \|x\|_2^2.$$

In \cite{connes-jones} a notion of property T was defined for
arbitrary von Neumann algebras (based on Kazhdan's property T for
groups).  It was observed that for i.c.c.\ groups $\Gamma$,
$L(\Gamma)$ has property T if and only if $\Gamma$ has Kazhdan's
property T.  In \cite[Theorem 3]{connes-jones} Connes and Jones also
prove the following theorem.

\begin{subthm}
\label{thm:T} No II$_1$-factor with property T can be embed into a
von Neumann algebra with Haagerup's approximation property.
\end{subthm}

\subsection{The Haagerup Invariant}

\begin{subdefn} If $M$ is a von Neumann algebra then we write
$\Lambda(M) \leq c$ if there exists a sequence of normal {\em finite
rank} c.b.\ maps $\Phi_k:M\to M$ such that (a) $\|\Phi_k\|_{cb} \leq
c$ for all $k$ and (b) $\xi(\Phi_k(x) - x) \to 0$ for all $x \in M$
and $\xi \in M_*$ (i.e.\ $\Phi_k \to id_M$ in the point $\sigma$-weak
topology). 

The {\em Haagerup Invariant} of $M$ is $$\Lambda(M) = \inf \{c:
\Lambda(M) \leq c \}.$$ Of course, $\Lambda(M) = \infty$ if no such $c$ 
exists.
\end{subdefn}

It is well known, and easily seen, that if $p \in M$ is a projection
then $\Lambda(pMp) \leq \Lambda(M)$ and if $M$ is a finite von Neumann
algebra and $N \subset M$ then $\Lambda(N) \leq \Lambda(M)$ (since
there is a conditional expectation $M\to N$).

In \cite{cowling-haagerup} Cowling and Haagerup showed, among other
things, that if $\Gamma$ is a lattice in $Sp(1,n)$ then
$\Lambda(L(\Gamma)) = 2n-1$ ($n\geq 2$).  Since $Sp(1,n)$ has
Kazhdan's property T it follows that $\Gamma$ also has property T and,
in particular, is finitely generated (cf.\ \cite{delaharpe-valette}).
A result of Malcev asserts that all finitely generated linear groups
are residually finite (cf. \cite{alperin}) and it is easy to see that
von Neumann algebras coming from residually finite groups always embed
into $R^{\omega}$.  Combining all of these facts we arrive at the
following theorem.

\begin{subthm} 
\label{thm:counterexample}
There is a finitely generated II$_1$-factor $N$ such
that $\Lambda(N) > 1$ and $N \subset R^{\omega}$.
\end{subthm}

Note that getting subfactors of $R^{\omega}$ just with property T is
much simpler than finding ones with non-trivial Haagerup invariant
since many property T groups (e.g.\ $SL(3,{\mathbb Z})$) are easily
seen to be residually finite.

We will also need the following theorem of Haagerup. 

\begin{subthm}(cf.\ \cite{haagerup:freegroup})
\label{thm:freegroup}
$\Lambda(L({\mathbb F}_n)) = 1$ for any free group ${\mathbb F}_n$ 
(on any number of generators).
\end{subthm}

\subsection{Weak Exactness for von Neumann Algebras}

In \cite{kirchberg} Kirchberg defined a notion of exactness in the
setting of von Neumann algebras.  He also listed a few facts without
proof.  In \cite{ozawa:weaklyexact} Ozawa gives proofs of these facts
as well as a few additional results.  Though it is actually a theorem
(cf.\ \cite{ozawa:weaklyexact}) we will define a von
Neumann algebra $M$ to be {\em weakly exact} if for each finite
dimensional operator system $X \subset M$ there exists a sequence of
u.c.p.\ maps $\phi_k:X\to M_{n(k)}({\mathbb C})$, $\psi_k:\phi_k(X)
\to M$ (note the domain of $\psi_k$!) such that
$$\xi(\psi_k\circ\phi_k(x) - x) \to 0,$$ for all $x \in X$ and all
$\xi \in M_*$.  Of course the proper definition is in terms of tensor
products but the theorem is that the finite dimensional approximation
property above is equivalent to the tensor product definition.  In any
case, we only need the approximation property. 

\begin{subthm}(cf.\ \cite{ozawa:weaklyexact}) 
\label{thm:exact}
If $\Gamma$ is discrete, $L(\Gamma)$ is weakly exact if and only
if $\Gamma$ is exact.
\end{subthm}
 
Neither direction of this theorem is easy but we wish to point out
that the ``if'' statement is very deep as it relies on Kirchberg's
remarkable theorem ``exactness implies local reflexivity''.

\section{Limes Inferior of von Neumann Algebras}

We recall a definition from \cite{haagerup-winslow:I}. 

\begin{defn} 
Let $M_k \subset B(H)$ be a sequence of von Neumann subalgebras of
$B(H)$. Then
$$\liminf_{k\to \infty} M_k = \{ x \in B(H): \exists x_k \in M_k \ 
\mathrm{such \ that} \ x_k \stackrel{so-*}{\to} x\},$$ where $x_k
\stackrel{so-*}{\to} x$ means convergence in the strong-$*$ topology.
\end{defn}

$\liminf_{k\to \infty} M_k$ is again a von Neumann algebra and in this
section we will observe that, in the special case that each $M_k$ is
contained in a fixed II$_1$-factor $M$, certain approximation
properties of the $M_k$'s will pass to $\liminf_{k\to \infty} M_k$. 

The following lemma isolates the main point. It is a simple matter to
extend the results of this section to  finite von Neumann
algebras but our interest is only in the factor case.

\begin{lem} Assume $M \subset B(L^2(M))$ is a II$_1$-factor in standard
form.  Let $M_k \subset M$ be von Neumann subalgebras and $E_k:M\to
M_k$ be the unique, normal, trace preserving conditional expectations
onto the $M_k$'s.  Then for each $x \in \liminf_{k\to \infty} M_k$,
$E_k(x) \stackrel{so-*}{\to} x$.
\end{lem}

\begin{proof} Since $E_k:M\to M_k \subset M$ extends to a contraction on
$L^2(M)$ we have $$\|x - E_k(x) \|_2 \leq \|x - x_k \|_2 +
\|E_k(x_k) - E_k(x) \|_2 \leq 2 \|x - x_k \|_2,$$ whenever $x_k \in
M_k$. Since $x_k \stackrel{so-*}{\to} x$ implies that $\|x - x_k \|_2
\to 0$ it follows that $\|x - E_k(x) \|_2 \to 0$ as well.  Since
$\{E_k(x) \}$ is norm bounded this further implies that $E_k(x) \to x$
in the strong operator topology.  Evidently the same argument shows
convergence of the adjoints too.
\end{proof}

\begin{rem} The lemma above and the proposition below are inspired by
the following exact C$^*$-algebra fact: If $A, A_k \subset B(H)$ are
C$^*$-algebras, each $a \in A$ is a norm limit of some sequence from
the $A_k$'s and each $A_k$ is exact then $A$ is also exact.  The proof
in the C$^*$-case uses Kirchberg's nuclear embeddability
characterization together with Arveson's Extension Theorem.  The point
of the lemma above is that in the II$_1$-factor context conditional
expectations can be used instead of Arveson's Theorem.
\end{rem}

\begin{prop} Assume $M \subset B(L^2(M))$ is a II$_1$-factor in standard
form and  $M_k \subset M$ are von Neumann subalgebras. 
\begin{enumerate}
\item If each $M_k$ is hyperfinite then so is $\liminf_{k\to \infty} M_k$.

\item If each $M_k$ is weakly exact then so is $\liminf_{k\to \infty} M_k$.

\item If each $M_k$ has the Haagerup approximation property then so 
does $\liminf_{k\to \infty} M_k$.

\item $\Lambda(\liminf_{k\to \infty} M_k) \leq \liminf_{k\to \infty} 
\Lambda(M_k)$.
\end{enumerate}
\end{prop}

\begin{proof} 
We only give the proofs of (2), (3) and (4) as these are the ones we will
use.  Moreover, the proof of (1) is similar to (and slightly easier
than) the proof of (2) using the equivalence of hyperfiniteness and
semi-discreteness (cf.\ \cite{connes:classification}). 

So assume that each $M_k$ is weakly exact.  We must show that for each
finite dimensional operator system $X \subset N = \liminf_{k\to
\infty} M_k$, finite set ${\mathfrak F} \subset X$, finite set
${\mathcal S} \subset N_*$ and each $\epsilon > 0$ there exist u.c.p.\
maps $\phi:X \to M_n({\mathbb C})$ and $\psi:\phi(X) \to N$ such that
$$|\xi(\psi\circ\phi(x)) - \xi(x) | < \epsilon$$ for all $x \in
{\mathfrak F}$ and $\xi \in {\mathcal S}$.  By the previous lemma we
can find a $k$ large enough that $$|\xi\circ E_N(E_k(x) - x)| <
\epsilon/2,$$ for all $x \in {\mathfrak F}$ and $\xi \in {\mathcal
S}$, where $E_N:M\to N$ is the trace preserving conditional
expectation.  Having fixed $k$ we can then find u.c.p.\ maps
$\phi_k:E_k(X) \to M_n({\mathbb C})$ and $\psi_k:\phi_k(E_k(X)) \to M_k$
such that $$|\xi\circ E_N\big(\psi_k\circ\phi_k(E_k(x)) - E_k(x)\big)| <
\epsilon/2,$$ for all $x \in {\mathfrak F}$ and $\xi \in {\mathcal S}$.  
Define $\phi:X\to M_n({\mathbb C})$ by $\phi = \phi_k\circ E_k|_X$
and $\psi:\phi(X) \to N$ by $\psi = E_N\circ \psi_k$. Then for each $x
\in {\mathfrak F}$ and $\xi \in {\mathcal S}$ we compute
\begin{eqnarray*}
|\xi(\psi\circ\phi(x)) - \xi(x) | 
&=&
|\xi\circ E_N \big( \psi_k\circ\phi_k(E_k(x)) - x \big)|\\
&\leq&
|\xi\circ E_N \big( \psi_k\circ\phi_k(E_k(x)) - E_k(x) \big)| + 
|\xi\circ E_N \big( E_k(x) - x \big)|\\ 
&<& 
\epsilon.
\end{eqnarray*}

To prove (3) we first fix a finite set ${\mathfrak F} \subset N =
\liminf_{k\to \infty} M_k$ and $\epsilon > 0$.  Choose $k$ large
enough that $\| E_k(x) - x \|_2$ is small for all $x \in {\mathfrak
F}$.  Then take a trace preserving u.c.p.\ map $\phi_k:M_k \to M_k$
such that $\| \phi_k(E_k(x)) - E_k(x) \|_2$ is small and $\phi_k$
induces compact operator on $L^2(M_k)$.  Evidently the composition
$E_N \circ \phi_k \circ E_k|_N$ induces a compact operator on $L^2(N)$
and one readily checks that $\| x - E_N \circ \phi_k \circ E_k(x)\|_2$
is small for all $x \in {\mathfrak F}$.

To prove the final assertion we first remark that it is clear from the
definition of limes inferior that if $\{k_j\}$ is any subsequence then
$$\liminf_{k\to \infty} M_k \subset \liminf_{j\to \infty} M_{k_j}.$$
It follows that $$\Lambda(\liminf_{k\to \infty}) \leq
\Lambda(\liminf_{j\to \infty} M_{k_j}),$$ and hence we may assume
$lim_{k\to \infty} \Lambda(M_k)$ exists and equals some number $c$.
With this observation the proof is quite similar to the arguments
above. Given finite sets ${\mathfrak F} \subset X$, ${\mathcal S}
\subset N_*$ and each $\epsilon > 0$ one chooses $k$ so that
$|\xi\circ E_N (E_k(x) - x)|$ is small (and $\Lambda(M_k) \leq c +
\epsilon$) then chooses finite rank, normal $\Phi_k:M_k \to M_k$ so
that $|\xi\circ E_N(\Phi_k(E_k(x)) - E_k(x))|$ is small (and
$\|\Phi_k\|_{cb} \leq c + \epsilon$) and, finally, defines $\Psi:N \to
N$ by $\Psi = E_N\circ \Phi_k \circ E_k$.
\end{proof}

\begin{rem} 
Since there exist subfactors $N \subset R^{\omega}$ with property T
and Haagerup invariant strictly larger than one (cf.\ Theorem
\ref{thm:counterexample}) it follows from \cite{haagerup-winslow:II}
that the analogues of statements (1) and (3) and (4)in the proposition
above do not hold when $M$ is replaced by $B(H)$. We don't know if the
analogue of (2) holds.
\end{rem}

\begin{cor} 
\label{thm:maincor}
Let $M \subset B(L^2(M))$ be a II$_1$-factor in standard form. Assume
$N, M_k \subset M$ are von Neumann subalgebras and further assume that
there exists a weakly dense $*$-subalgebra $A \subset N$ such that $A
\subset \liminf_{k\to \infty} M_k$. Then
$$\Lambda(N) \leq \liminf_{k\to \infty} \Lambda(M_k).$$ If
each $M_k$ is weakly exact or has the Haagerup
approximation property then so does $N$.
\end{cor}

\begin{proof} Since $\liminf_{k\to \infty} M_k$ is a von Neumann
algebra it follows that $N \subset \liminf_{k\to \infty} M_k$.  The
result then follows from the previous proposition together with the
fact that in finite von Neumann algebras weak exactness or the
Haagerup approximation property passes to subalgebras and Haagerup
invariants do not increase.
\end{proof}

\section{Finite Free Entropy}

\begin{thm} 
\label{thm:mainthm}
There exist noncommutative random variables $X_1,\ldots, X_n$ with the
property that $\chi(X_1,\ldots, X_n) > -\infty$ but $M =
W^*(X_1,\ldots, X_n)$ is not isomorphic to any (not necessarily
unital) subalgebra of a free group factor.
\end{thm}

\begin{proof} 
Let $N \subset R^{\omega}$ be a finitely generated II$_1$-factor with
Haagerup invariant strictly larger than one or just a property T
factor (cf.\ Theorem \ref{thm:counterexample}).  Let $Y_1,\ldots,Y_n$
be a set of self-adjoint generators of $N$.  Inside $N\ast L({\mathbb
F}_n)$ consider the self-adjoints $$Y_1 + \epsilon S_1,\ldots,Y_n +
\epsilon S_n,$$ where $S_i \in L({\mathbb F}_n)$ are free semicircular
elements, and let $N_{\epsilon}$ be the von Neumann algebra generated
by these elements.

By the free entropy-power inequality \cite[Theorem
3.9]{voiculescu:asymptotic} we have that $$\chi(Y_1 + \epsilon
S_1,\ldots,Y_n + \epsilon S_n) > -\infty$$ for every $\epsilon >
0$. (Strictly speaking this also uses the fact that free elements are
regular in the sense of \cite[Definition 3.6]{voiculescu:asymptotic}.)
Now it is clear that every element in the $*$-algebra generated by
$Y_1,\ldots,Y_n$ is a norm limit of elements from $N_{\epsilon}$ and
hence we can apply Corollary \ref{thm:maincor} to conclude that
$$\liminf \Lambda(N_{\epsilon}) > 1.$$ Hence for all sufficiently
small $\epsilon$ it follows that $N_{\epsilon}$ cannot be embed (even
non-unitally) into a free group factor (cf.\ Theorem
\ref{thm:freegroup}).  Of course, if one just started with a property T
factor then the same conclusion holds since all free group factors
have the Haagerup approximation property.
\end{proof}

\begin{rem} Ozawa has pointed out that if one starts with a co-compact
lattice $\Gamma \subset Sp(1,n)$ then $\Gamma$ will be hyperbolic.
If follows that $\Gamma \ast {\mathbb F}_n$ is also hyperbolic and,
since subfactors of solid factors are again solid, one can arrange the
example above to be solid \cite{ozawa:solid}.
\end{rem}

In \cite{kenley} a notion of free Hausdorff dimension ${\mathbb
H}(\cdot)$ was introduced and it was observed that ${\mathbb
H}(X_1,\ldots, X_n) \leq \delta_0(X_1,\ldots, X_n)$, where $\delta_0$
is the ``modified'' free entropy dimension (cf.\ \cite[Corollary
3.5]{kenley}).  In particular one has the following general inequality
$${\mathbb H}(X_1,\ldots, X_n) \leq \delta(X_1,\ldots, X_n).$$ By the
results of \cite{bcg} we have that $$\delta(X_1,\ldots, X_n) \leq
\delta^*(X_1,\ldots, X_n),$$ where $\delta^*$ is Voiculescu's
``non-microstates'' free entropy dimension
\cite{voiculescu:survey}. Finally, in \cite{connes-dima} a numerical
invariant $\Delta(\cdot)$ was introduced and it was shown that (cf.\
\cite[Theorems 3.3 and 4.4 and Lemma 4.1]{connes-dima})
$$\delta^*(X_1,\ldots, X_n) \leq \Delta(X_1,\ldots, X_n) \leq n.$$
Putting all of these estimates together we have $${\mathbb
H}(X_1,\ldots, X_n) \leq \delta(X_1,\ldots, X_n) \leq
\delta^*(X_1,\ldots, X_n) \leq \Delta(X_1,\ldots, X_n) \leq n$$ for
arbitrary random variables in a tracial W$^*$-probability space.

\begin{cor} 
There exist noncommutative random variables $X_1,\ldots, X_n$ with the
property that $${\mathbb H}(X_1,\ldots, X_n) = \delta(X_1,\ldots, X_n)
= \delta^*(X_1,\ldots, X_n) = \Delta(X_1,\ldots, X_n) = n$$ 
but $M = W^*(X_1,\ldots, X_n)$ is not isomorphic to a subalgebra of a 
free group factor.
\end{cor}

\begin{proof} \cite[Corollary 3.8]{kenley} states that if
$\chi(X_1,\ldots, X_n) > -\infty$ then ${\mathbb H}(X_1,\ldots, X_n) =
n$ and hence the examples from the previous theorem satisfy the
conclusion of  this corollary.
\end{proof}

\section{Free Probability and the Novikov Conjecture?}

Consider a much weaker version of the question
addressed in the previous section.  
\begin{center}
{\em If $\chi(X_1,\ldots,X_n) > -\infty$ is 
$W^*(X_1,\ldots,X_n)$ a weakly exact von Neumann algebra?}
\end{center}

Note that this is significantly weaker than the problem answered
(negatively) in this paper.  Indeed, free group factors are weakly
exact and at present we do not have any {\em explicit} examples of
non-weakly exact von Neumann algebras. (cf.\ Theorem \ref{thm:exact}
-- Gromov's non-exact groups yield non-weakly exact von Neumann
algebras but Gromov's examples arise as inductive limits and are not
very ``concrete'' \cite{ollivier}).

In \cite[Question 8, Section 12]{brown:invariantmeans} we posed the
following problem: If $L(\Gamma) \subset R^{\omega}$ does it follow
that $\Gamma$ is exact or, perhaps, just uniformly embeddable into
Hilbert space?  An affirmative answer to this question would have
three important consequences.  Namely, it would imply (1) not every
II$_1$-factor embeds into $R^{\omega}$ (hence a negative answer to
Voiculescu's ``unification problem''), (2) not every hyperbolic group
is residually finite (a longstanding problem in geometric group
theory) and (3) the Novikov conjecture holds for all residually finite
(or $R^{\omega}$-embeddable) groups. (See \cite{brown:invariantmeans} 
for a more detailed discussion.)

\begin{prop} If ``$\chi(X_1,\ldots,X_n) > -\infty \Longrightarrow
W^*(X_1,\ldots,X_n)$ is weakly exact'' then every discrete group which embeds
into $R^{\omega}$ is exact.
\end{prop}

It is known that if a group admits a faithful unitary representation
into $R^{\omega}$ then, in fact, $L(\Gamma) \subset R^{\omega}$ (cf.\
\cite{kirchberg:T},\cite{radulescu}).  Hence, assuming the implication
``$\chi(X_1,\ldots,X_n) > -\infty \Longrightarrow W^*(X_1,\ldots,X_n)$
is weakly exact'', we could approximate $L(\Gamma) \subset L(\Gamma)
\ast L({\mathbb F}_n)$ by weakly exact von Neumann algebras just as we
did in the proof of Theorem \ref{thm:mainthm}.  It would follow from
Corollary \ref{thm:maincor} that $L(\Gamma)$ was also weakly exact and
hence, by Theorem \ref{thm:exact}, that $\Gamma$ is exact.  In case
one worries about finite generation we remark that a discrete group is
exact if and only if all of its finitely generated subgroups are
exact.

We should mention that it is possible that finite free entropy may not
imply any sort of finite dimensional approximation property (other
than the obvious one of large microstate spaces).  More precisely, if
it turns out that every hyperbolic group embeds into $R^{\omega}$ then
it would follow that Gromov's non-exact groups also embed into
$R^{\omega}$ and hence there would exist a non-weakly exact
II$_1$-factor $M = W^*(X_1,\ldots,X_n)$ such that
$\delta(X_1,\ldots,X_n) = n$.

\bibliographystyle{amsplain}

\providecommand{\bysame}{\leavevmode\hbox to3em{\hrulefill}\thinspace}

\end{document}